\documentclass[12pt]{article}

\usepackage{e-jc}
\usepackage{latexsym,amsthm,amssymb,amsmath}
\usepackage{epsfig,graphicx,color}

\newtheorem{theorem}{Theorem}[section]
\newtheorem{corollary}[theorem]{Corollary}

\title{The spectral radius and the maximum degree of irregular graphs}

\author{Sebastian M. Cioab\u{a}\footnote{Research partially supported by an NSERC postdoctoral fellowship.}
\\
\small Department of Mathematics\\[-0.8ex]
\small University of California, San Diego\\[-0.8ex]
\small La Jolla, CA 92093-0112\\[-0.8ex]
\small \texttt{scioaba@math.ucsd.edu}}

\date{\dateline{Jan 20, 2007}{??}\\
\small MR Subject Classifications: 05C50, 15A18}

\begin{document}
\maketitle

\begin{abstract}
Let $G$ be an irregular graph on $n$ vertices with maximum degree
$\Delta$ and diameter $D$. We show that
\begin{equation*}
\Delta-\lambda_1>\frac{1}{nD}
\end{equation*}
where $\lambda_1$ is the largest eigenvalue of the adjacency matrix
of $G$. We also study the effect of adding or removing few edges on
the spectral radius of a regular graph.
\end{abstract}

\section{Preliminaries}

Our graph notation is standard (see West \cite{W}). For a graph $G$,
we denote by $\lambda_i(G)$ the $i$-th largest eigenvalue of its
adjacency matrix and we call $\lambda_1(G)$ the spectral radius of
$G$. If $G$ is connected, then the positive eigenvector of norm $1$
corresponding to $\lambda_1(G)$ is called the principal eigenvector
of $G$.

The spectral radius of a connected graph has been well studied.
Results in the literature connect it with the chromatic number, the
independence number and the clique number of a connected graph
\cite{GN,Hae1,Hoff1,Nik2,Wilf1}. Recently, it has been shown that
the spectral radius also plays an important role in modeling virus
propagation in networks \cite{JKMD,WCWF}.

In this paper, we are interested in the connection between the
spectral radius and the maximum degree $\Delta$ of a connected graph
$G$. In particular, we study the spectral radius of graphs obtained
from $\Delta$-regular graphs on $n$ vertices by deleting a small
number of edges or loops. The Erd\H{o}s-R\'{e}nyi graph $ER(q)$ is
an example of such a graph, see \cite{GN,N} and the references
within for more details on its spectral radius and other interesting
properties.

It is a well known fact that $\lambda_1(G)\leq \Delta(G)$ with
equality if and only if $G$ is regular. It is natural to ask how
small $\Delta(G)-\lambda_1(G)$ can be when $G$ is irregular.

Cioab\u{a}, Gregory and Nikiforov \cite{CGN} proved that if $G$ is
an irregular graph on $n$ vertices, with maximum degree $\Delta$ and
diameter $D$, then
\begin{equation*}
\Delta-\lambda_1>\frac{1}{n(D+\frac{1}{n\Delta-2m})}\geq
\frac{1}{nD+n}
\end{equation*}
where $m$ is the number of edges of $G$. This result improved
previous work of Stevanovi\'{c} \cite{Stev}, Zhang \cite{Z} and Alon
and Sudakov \cite{AS}.

In \cite{CGN}, the authors conjecture that
\begin{equation}\label{cgn_conj}
\Delta-\lambda_1>\frac{1}{nD}
\end{equation}

In this paper, we prove this conjecture. Using inequality
\eqref{cgn_conj}, we improve some recent results of Nikiforov
\cite{Nik1} regarding the spectral radius of a subgraph of a regular
graph. We also investigate the spectral radius of a graph obtained
from a regular graph by adding an edge.

For recent results connecting the spectral radius of a general (not
necessarily regular) graph $G$ and that of a subgraph of $G$, see
\cite{Nik1}.

\section{The spectral radius and the maximum degree}

The following theorem is the main result of this section.
\begin{theorem}\label{nD}
Let $G$ be a connected irregular graph with $n$ vertices, maximum
degree $\Delta$ and diameter $D$. Then
\begin{equation*}
\Delta-\lambda_1(G)>\frac{1}{nD}
\end{equation*}
\end{theorem}
\begin{proof}
Let $x$ be the principal eigenvector for $G$. Let $s$ be a vertex of
$G$ such that $x_s=\max_{i\in [n]}x_i$. Since $G$ is not regular, it
follows that $x_s>\frac{1}{\sqrt{n}}$.

If the degree of $s$ is not $\Delta$, then
$$
\lambda_1 x_s=\sum_{j\sim s}x_j\leq (\Delta-1)x_s
$$
which implies $\Delta-\lambda_1\geq 1>\frac{1}{nD}$ and proves the
theorem.

From now on, we will assume that the degree of $s$ is $\Delta$.

Suppose first that $G$ contains at least two vertices whose degree
is not $\Delta$. Let $u$ and $v$ be two vertices of $G$ whose degree
is not $\Delta$.

Let $P:u=i_0,i_1,\dots, i_r=s$ be a shortest path from $u$ to $s$ in
$H$. Obviously, $r\leq D$. Let $Q$ be a shortest path from $v$ to
$s$ in $H$. Let $t$ be the smallest index $j$ such that $i_j$ is on
$Q$. Obviously, $t\in \{0,\dots, r\}$.

If $t=0$, then the distance from $u$ to $s$ is at most $D-1$ (this
means $r\leq D-1$), and applying a similar argument to the one in
\cite{CGN}, we obtain that
\begin{align*}
\Delta-\lambda_1(G)&=\Delta \sum_{j=1}^{n}x_{j}^{2}-\sum_{kl\in
E(G)}2x_kx_l\\
&=\sum_{i=1}^{n}(\Delta-d_i)x_i^{2}+\sum_{kl\in E(G)}(x_k-x_l)^{2}\\
&\geq x_u^{2}+\sum_{j=0}^{r-1}(x_{i_{j+1}}-x_{i_j})^{2}\\
&\geq
\frac{\left(x_u+\sum_{j=0}^{r-1}(x_{i_{j+1}}-x_{i_j})\right)^{2}}{r+1}=\frac{x_{i_r}^{2}}{r+1}>\frac{1}{nD}
\end{align*}

If $t\geq 1$, we may assume without any loss of generality that
$t=d(u,i_t)\geq d(v,i_t)$. Let $Q_{v,i_t}$ denote the sub-path of
$Q$ which connects $v$ to $i_t$. Using the Cauchy-Schwarz
inequality, it follows that
\begin{align*}
\Delta-\lambda_1(G)&=\sum_{i=1}^{n}(\Delta-d_i)x_i^{2}+\sum_{kl\in E(G)}(x_k-x_l)^{2}\geq x_u^{2}+x_v^{2}+\sum_{kl\in E(G)}(x_k-x_l)^{2}\\
&\geq
(x_u^{2}+\sum_{j=0}^{t-1}(x_{i_j}-x_{i_{j+1}})^{2})+(x_v^{2}+\sum_{kl\in
E(Q_{v,i_t})}(x_k-x_l)^{2})+\sum_{j=t}^{r}(x_{i_j}-x_{i_{j+1}})^{2}\\
&\geq
\frac{x_{i_t}^{2}}{t+1}+\frac{x_{i_{t}}^{2}}{d(v,i_t)+1}+\frac{(x_{i_t}-x_s)^{2}}{r-t}\\
&\geq \frac{2x_{i_t}^{2}}{t+1}+\frac{(x_s-x_{i_t})^{2}}{r-t}
\end{align*}
The right hand-side is a quadratic function in $x_{i_t}$ which
attains its minimum when $x_{i_t}=\frac{(t+1)x_{s}}{2r-t+1}$. This
implies that
\begin{equation*}
\Delta-\lambda_1(G)>\frac{2x_s^{2}}{2r-t+1}\geq \frac{x_s^{2}}{r}
\end{equation*}
since $t\geq 1$. Because $x_s>\frac{1}{\sqrt{n}}$ and $r\leq D$, we
obtain
\begin{equation*}
\Delta-\lambda_1(G)>\frac{1}{nD}
\end{equation*}
This finishes the proof in the case that $G$ has at least two
vertices whose degree is not $\Delta$.

Assume now that $G$ contains exactly one vertex whose degree is less
than $\Delta$. Let $w$ be a vertex whose principal eigenvector entry
is minimum. Then $d_{w}<\Delta$ because
$$
\Delta x_{w}>\lambda_1 x_{w}=\sum_{j\sim w}x_{j}\geq d_{w}x_{w}
$$

Recall that $x_{s}=\max_{i\in [n]}x_{i}$. Let
$\gamma=\frac{x_{s}}{x_{w}}$. We may assume that $\gamma>D$.
Otherwise, by summing the equalities $\lambda_1 x_{i}=\sum_{j\sim
i}x_{j}$ over all $i\in [n]$ we have
$$
\Delta-\lambda_1=\frac{(\Delta-d_{w})x_{w}}{\sum_{i=1}^{n}x_{i}}>\frac{x_{w}}{nx_{s}}=
\frac{1}{n\gamma}\geq \frac{1}{nD}
$$
which proves the theorem.

We may also assume that $d(w,s)=D$ because otherwise by applying an
argument similar to the one of the previous case, we can prove the
theorem.

We claim there exists $j\sim s$ such that $x_j<\frac{1}{\sqrt{n}}$.
Otherwise, let $j\sim s$ such that $d(j,w)=D-1$. Then applying the
argument from the previous case gives
$$
\Delta-\lambda_1 >\frac{x_{j}^{2}}{D}>\frac{1}{nD}
$$
which again proves the theorem.

Since $x_j<\frac{1}{\sqrt{n}}$ and $j\sim s$, we have
$$
\lambda_1 x_{s}=\sum_{l\sim
s}x_{l}<(\Delta-1)x_{s}+\frac{1}{\sqrt{n}}
$$
which implies
$$
\Delta-\lambda_1>1-\frac{1}{x_{s}\sqrt{n}}
$$
If the right-hand side is at least $\frac{1}{nD}$, then we are done.
Otherwise, $1-\frac{1}{x_{s}\sqrt{n}}<\frac{1}{nD}$ implies
\begin{equation}\label{x_s}
x_{s}<\frac{D\sqrt{n}}{nD-1}
\end{equation}
Since $x_{s}=\max_{i\in [n]}x_{i}, x_{w}=\min_{i\in [n]}x_{i}$, we
have that
$$
(n-1)x_{s}^{2}+x_{w}^{2}\geq \sum_{l=1}^{n}x_{l}^{2}=1
$$
which implies
\begin{equation}\label{x_w}
x_{w}^{2}\geq 1-\frac{(n-1)nD^{2}}{(nD-1)^{2}}=\frac{
(nD-1)^{2}-(n-1)nD^{2}}{(nD-1)^{2}}=\frac{(D^{2}-2D)n+1}{(nD-1)^{2}}
\end{equation}
Assume $D\geq 3$. From \eqref{x_s} and \eqref{x_w}, we get that
$$
\gamma^{2}=\frac{x_{s}^{2}}{x_{w}^{2}}<\frac{D^{2}n}{(D^{2}-2D)n+1}<D^{2}
$$
Thus, $\gamma<D$ which is a contradiction with the earlier
assumption that $\gamma>D$. This proves the theorem for $D\geq 3$.

For $D=2$, looking at the square of the adjacency matrix of $G$ we
get
\begin{equation}\label{square}
\lambda_1^{2} x_{s}\leq (\Delta^{2}-1)x_{s}+ x_{w}
\end{equation}
which implies
$$
\lambda_1^{2}\leq \Delta^{2}-1+\frac{1}{\gamma}<
\Delta^2-\frac{1}{2}
$$
since $\gamma\geq 2$. Note that inequality \eqref{square} holds
because there is at least one path of length $2$ from $s$ to $w$.

Thus, $\lambda_1\leq
\sqrt{\Delta^{2}-\frac{1}{2}}<\Delta-\frac{1}{4\Delta}$

If $n\geq 2\Delta$, then we are done. Suppose then that $n<2\Delta$.
Then the vertex $s$ has at least two neighbours at distance $1$ from
the vertex $w$. We deduce that
$$
\lambda_1^{2}x_{s}\leq (\Delta^{2}-2)x_{s}+2x_{w}
$$
which implies
$$
\lambda_1^{2}\leq \Delta^{2}-2+\frac{2}{\gamma}\leq \Delta^{2}-1
$$
Thus,
$$
\lambda_1\leq
\sqrt{\Delta^{2}-1}<\Delta-\frac{1}{2\Delta}<\Delta-\frac{1}{nD}
$$
which completes the proof of the theorem.
\end{proof}

Because $\lambda_1(G)=\Delta(G)$ when $G$ is regular, the following
result is an immediate consequence of Theorem \ref{nD}.
\begin{corollary}
Let $G$ be a $\Delta$-regular graph and $e$ be an edge of $G$ such
that $G\setminus e$ is connected. Then
\begin{equation*}
\frac{2}{n}>\Delta-\lambda_1(G\setminus e)>\frac{1}{nD}
\end{equation*}
where $D$ is the diameter of $G\setminus e$.
\end{corollary}

The previous results improve Theorems 4-6 obtained by Nikiforov in a
recent paper \cite{Nik1}.

If $f,g:\mathbb{N}\rightarrow [0,+\infty)$ we write $f(n)=O(g(n))$
if there is $c>0$ and $n_0>0$ such that $f(n)\leq cg(n)$ for $n\geq
n_0$ and we write $f(n)=\Theta(g(n))$ if $f(n)=O(g(n))$ and
$g(n)=O(f(n))$.

Under the same hypothesis as the previous corollary, if $\Delta$ is
fixed and $G\setminus e$ is connected, then the diameter $D$ of
$G\setminus e$ is at least $\log_{\Delta-1}n+O(1)$. In this case, we
obtain the following estimates
$$
\frac{2}{n}>\Delta-\lambda_1(G\setminus
e)>O\left(\frac{1}{n\log_{\Delta-1}n}\right)
$$
It seems likely that the upper bound gives the right order of
magnitude for $\Delta-\lambda_1(G\setminus e)$, but proving this
fact is an open problem.

If $\Delta=\Theta(n)$ and $G\setminus e$ is connected, then the
diameter $D$ of $G\setminus e$ is $O(1)$. To see this consider a
path $i_0,\dots, i_D$ of length $D$ in $G\setminus e$. For $0\leq
j\leq D$, let $N_j$ denote the neighborhood of vertex $j$ in
$G\setminus e$. It follows that for $j\equiv 0\pmod{3}, 0\leq j\leq
D$, the sets $N_{j}$ are pairwise disjoint. Also, for all but at
most two $j$'s, we have $|N_j|=\Delta$. These facts imply that
$$
n\geq \sum_{j\equiv
0\pmod{3}}|N_j|>\left(\frac{D+1}{3}-2\right)\Delta+2(\Delta-1)=\frac{\Delta(D+1)}{3}-2
$$
Thus,
$$
D<\frac{3(n+2)}{\Delta}=O(1)
$$
Hence, in this case, our estimates imply that
$$
\Delta-\lambda_1(G\setminus e)=\Theta\left(\frac{1}{n}\right)
$$

Note that the previous argument can be also used to show that if
$\Delta=\Theta(n)$ and $H$ is a connected graph obtained from a
$\Delta$-regular graph on $n$ vertices by deleting a constant number
of edges, then
$$
\Delta-\lambda_1(H)=\Theta\left(\frac{1}{n}\right)
$$

\section{Adding an edge to a regular graph}

In this section, we analyze the effect of adding an edge on the
spectral radius of a regular graph. We need different techniques in
this case because the spectral radius will be closer to the minimum
degree than to the maximum degree of the graph.
\begin{theorem}
Let $H$ be a connected, $k$-regular graph and $e\notin E(H)$. If $G$
is the graph obtained from $H$ by adding the edge $e$ and
$k-\lambda_2(H)>1$, then
\begin{equation}
\frac{2}{n}\cdot
\left(1+\frac{1}{k-\lambda_2(H)-1}\right)>\lambda_1(G)-\lambda_1(H)>\frac{2}{n}\cdot
\left(1+\frac{1}{2(k+1)}\right)
\end{equation}
\end{theorem}
\begin{proof}
The lower bound follows by applying the following result obtained by
Nikiforov \cite{Nik3}. See also \cite{CG} for related results.
\begin{theorem}
Let $G$ be an irregular graph with $n$ vertices and $m$ edges having
maximum degree $\Delta$. If $G$ has at least two vertices of degree
$\Delta$ and at least two vertices of degree less than $\Delta$,
then
$$
\lambda_1(G)>\frac{2m}{n}+\frac{2}{4m+1}
$$
\end{theorem}

For the upper bound, we use the following result of Maas
\cite{Maas1} (see also Theorem 6.4.1 in \cite{CRS} and \cite{ZB}).
\begin{theorem}[Maas \cite{Maas1}]
Let $x$ be the principal eigenvector of a graph $H$ and let $i$ and
$j$ be two non-adjacent vertices in $H$. Then
\begin{equation}
\lambda_1(H+ij)-\lambda_1(H)<1+\delta-\beta
\end{equation}
where
$$
\beta=\lambda_1(H)-\lambda_2(H)
$$
and $\delta$ satisfies the equation
$$
\frac{\delta(1+\delta)(2+\delta)}{(x_i+x_j)^{2}+\delta(2+\delta+2x_ix_j)}=\beta
$$
\end{theorem}

By applying the previous theorem and using the fact that the
principal eigenvector of $H=G\setminus e$ has all entries equal to
$\frac{1}{\sqrt{n}}$, we obtain that
$$
\lambda_1(G)-\lambda_1(H))<1+\delta-\beta=\frac{2\beta}{\delta n}
$$

If $k-\lambda_2(H)=\beta >1$, we obtain that
$$
\beta=\frac{\delta(1+\delta)(2+\delta)}{(\delta+2)\left(\delta+\frac{2}{n}\right)}=\frac{\delta(1+\delta)}{\delta+\frac{2}{n}}<1+\delta
$$
Thus,
$\frac{2\beta}{\delta}<\frac{2\beta}{\beta-1}=\frac{2(k-\lambda_2(H))}{k-\lambda_2(H)-1}$.
Hence, we deduce that
$$
\lambda_1(G)-\lambda_1(G\setminus
e)<\frac{2(k-\lambda_2(H))}{(k-\lambda_2(H)-1)n}
$$
which proves the theorem.
\end{proof}

For $k\geq 3$ fixed and $\epsilon>0$, Friedman \cite{Fri} proved
that most $k$-regular graphs $H$ have $\lambda_2(H)\leq
2\sqrt{k-1}+\epsilon$. This implies that for most $k$-regular graphs
$H$ and for each $e\notin E(H)$, we have
\begin{equation}\label{add_edge}
\lambda_1(H+e)-\lambda_1(H)=\Theta\left(\frac{1}{n}\right)
\end{equation}

If $H$ is a $k$-regular graph with $k-\lambda_2(H)\leq 1$, then
\eqref{add_edge} might not hold. This is true at least for $k=2$ as
seen by the graph $G_n$ on $n$ vertices which is obtained from a
cycle on $n$ vertices by adding an edge between two vertices at
distance $2$. It follows from the work of Simi\'{c} and Koci\'{c}
\cite{SK} (see also \cite{CRS} equation (3.4.5) on page 63) that
$$
\lim_{n\rightarrow \infty}\lambda_1(G_n)=2.3829
$$
while $\lambda_1(C_n)=2$ so clearly
$\lambda_1(G_n)-\lambda_1(C_n)\neq \Theta\left(\frac{1}{n}\right)$.

\section{Final Remarks}

It is worth mentioning that there are infinite families of irregular
graphs with maximum degree $\Delta$ such that $\Delta-\lambda_1\leq
\frac{c}{nD}$, where $c$ is an absolute constant. Cioab\u{a},
Gregory and Nikiforov \cite{CGN} describe such a family with
$c=4\pi^{2}$ while Liu, Shen and Wang \cite{LSW} found an infinite
family with $c=3\pi^{2}$. It has yet to be determined what the best
such constant $c$ can be for all $n$ and $D$.

Note that the argument of Theorem \ref{nD} can be extended easily to
multigraphs. From the proof of Theorem \ref{nD}, it is easy to see
that if one can find better upper bounds for the distance between
vertices corresponding to the extreme entries of the principal
eigenvector of an irregular graph $G$, then one can improve the
result of Theorem \ref{nD}. However, there are infinite families of
irregular graphs for which the distance between such vertices equals
the diameter. We describe one such family in the next paragraph.

For $k\geq 2$, consider the cycle on $n=2k+1$ vertices with vertex
set $\{1,\dots,2k+1\}$ and edges $\{i,i+1\}$ for $1\leq i\leq 2k$
and $\{2k+1,1\}$. Add the edges $\{k,k+2\}, \{k+1,k+3\}$ and for
$k\geq 3$, the edges $\{i,2k+3-i\}$ for $2\leq i\leq k-1$ . The
resulting graph $G$ has maximum degree $3$ and $1$ is the only
vertex whose degree is $2$. It can be shown easily by induction that
if $x$ is the principal eigenvector of $G$, then $x_{1}=\min_{i\in
[n]}x_{i}$ and $x_{k+1}=x_{k+2}=\max_{i\in [n]}x_{i}$. Note that
$d(1,k+1)=d(1,k+2)=k$ which equals the diameter of $G$.

When $n=7$, the graph obtained by the above procedure is shown in
Figure 1 with its vertices labeled by their entries in the principal
eigenvector. Note that although the distance between vertices
corresponding to extreme entries of the principal eigenvector equals
the diameter, there are vertices whose eigenvector entry is at least
$\frac{1}{\sqrt{7}}=0.378$ which are at distance less than the
diameter from the vertex whose eigenvector entry is minimum. If this
fact would be true for any irregular graph, it would imply Theorem
\ref{nD}.

\begin{figure}[ht]
\begin{center}
\includegraphics{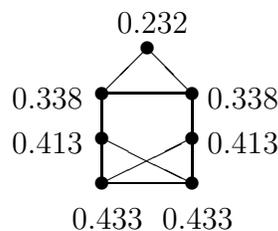}
\caption{An irregular graph and its principal eigenvector}
\end{center}
\end{figure}

It would be interesting to determine the precise asymptotic
behaviour of the spectral radius of a graph obtained from a
$\Delta$-regular graph on $n$ vertices by deleting an edge when
$\Delta=o(n)$. Another problem of interest is to find the exact
asymptotic behaviour of the spectral radius of a graph obtained from
a $k$-regular graph $H$ by adding an edge when $k=\Theta(n)$ or when
$k-\lambda_2(H)\leq 1$.

A slightly different direction of research was taken by Biyikoglu
and Leydold in \cite{BL} where they study the graphs which have the
maximum spectral radius in the set of all connected graphs with
given degree sequence. In particular, the authors show that the
maximum is increasing with respect to the majorization order. Even
for graphs with simple degree sequences, determining the maximum
spectral radius seems a nontrivial problem.

\section*{Acknowledgments} I am grateful to Steve Butler and David Gregory
for their careful reading of the paper and I thank Fan Chung, Orest
Bucicovschi and Vlado Nikiforov for helpful discussions.

\end{document}